%% file: socg.tex
\documentclass[a4paper,UKenglish,cleveref, autoref, thm-restate,nolineno]{socg-lipics-v2021}

\newboolean{appendix} 
\setboolean{appendix}{true} 

\title{An Improved
Lower Bound on the Number of Pseudoline Arrangements\footnote{%
This article is a
merge of the following two independent SoCG'24 submissions:\\
1) Justin Dallant
-- Improved Lower Bound on the Number of Pseudoline Arrangements
\cite{DallantVersion} \\
2) Fernando Cortés Kühnast, Stefan Felsner, and Manfred Scheucher --
An Improved Lower Bound on the Number of Pseudoline Arrangements 
\cite{CKFSVersion}.
First steps were made in the Bachelor's thesis of the first author \cite{FernandoBakkArbeit2023}.\\
While in 1) and 2) the estimates $b_n \ge 0.2604$ and $b_n \ge 0.2526$ were shown, 
we here establish $b_n \ge \cOptTwelve$ by combining our two approaches.
}} 

\author{Fernando {Cortés Kühnast}}{Institute of Mathematics, Technische Universität Berlin, Germany}{cortes@math.tu-berlin.de}{https://orcid.org/0009-0008-1847-362X}{}

\author{Justin Dallant
}{Algorithms Research Group, Université libre de Bruxelles, Belgium}{justindallant@gmail.com}{https://orcid.org/0000-0001-5539-9037}{}

\author{Stefan Felsner}{Institute of Mathematics, Technische Universität Berlin, Germany}{felsner@math.tu-berlin.de}{https://orcid.org/0000-0002-6150-1998}{}

\author{Manfred Scheucher}{Institute of Mathematics, Technische Universität Berlin, Germany}{scheucher@math.tu-berlin.de}{https://orcid.org/0000-0002-1657-9796}{}

\authorrunning{F. Cortés Kühnast, J. Dallant, S. Felsner, and M. Scheucher} 

\Copyright{Fernando Cortés Kühnast, Justin Dallant, Stefan Felsner, and Manfred Scheucher} 



\ifthenelse{\boolean{appendix}}
{}
{\relatedversion{A full version of this article is available at \url{https://arxiv.org/abs/2402.13107}}}

\supplement{Our supplemental data is available on github \cite{supplemental_data}}

\funding{
S.~Felsner was partially supported by DFG Grant FE~340/13-1. 
F.~Cortés Kühnast and M.~Scheucher were supported by DFG Grant SCHE~2214/1-1.
J.~Dallant was supported by the French Community of Belgium via the funding of a FRIA grant.}

\EventEditors{Wolfgang Mulzer and Jeff M. Phillips}
\EventNoEds{2}
\EventLongTitle{40th International Symposium on Computational Geometry (SoCG
2024)}
\EventShortTitle{SoCG 2024}
\EventAcronym{SoCG}
\EventYear{2024}
\EventDate{June 11-14, 2024}
\EventLocation{Athens, Greece}
\EventLogo{}
\SeriesVolume{293}
\ArticleNo{XX}     

\usepackage{enumerate}
\usepackage[style=english]{csquotes}

\graphicspath{{figs/}}

\usepackage[font=small,labelfont=bf]{caption}
\usepackage[font=small,labelfont=normalfont,labelformat=simple]{subcaption}
\renewcommand\thesubfigure{(\alph{subfigure})}

\usepackage{graphicx}

\usepackage{xcolor} 
\usepackage{todonotes}

\newcommand{\clow}{c^-}
\newcommand{\cupp}{c^+}

\DeclareMathOperator{\area}{area}

\newcommand{\floorfrac}[2]{{\left\lfloor\frac{#1}{#2} \right\rfloor}}
\newcommand{\crossing}[2]{{{#1} \times {#2}}}

\newcommand{\calL}{\mathcal{L}}

\newcommand{\R}{\mathbb{R}}
\newcommand{\cOptTwelve}{0.2721}

\newcommand{\cOptSix}{0.2542}
\newcommand{\cOptFour}{0.2183}

\definecolor{color1} {rgb}{0.133,0.545,0.133}
\definecolor{color2} {rgb}{0.000,0.000,0.545}
\definecolor{color3} {rgb}{0.690,0.188,0.376}
\definecolor{color4} {rgb}{1.000,0.271,0.000}
\definecolor{color5} {rgb}{1.000,1.000,0.000}
\definecolor{color6} {rgb}{0.871,0.722,0.529}
\definecolor{color7} {rgb}{0.000,1.000,0.000}
\definecolor{color8} {rgb}{0.000,1.000,1.000}
\definecolor{color9} {rgb}{1.000,0.000,1.000}
\definecolor{color10}{rgb}{0.392,0.584,0.929}

\begin{document}

\maketitle

\begin{abstract}
Arrangements of pseudolines are classic objects in discrete and computational geometry. They have been studied with increasing intensity since
their introduction almost 100 years ago.
The study of the number $B_n$ of non-isomorphic simple arrangements of $n$ pseudolines goes back to Goodman and Pollack, Knuth, and others. 
It is known that $B_n$ is in the order of $2^{\Theta(n^2)}$
and finding asymptotic bounds on $b_n = \frac{\log_2(B_n)}{n^2}$ remains a challenging task.
In 2011, Felsner and Valtr 
showed that $0.1887 \leq b_n \le 0.6571$ for sufficiently large~$n$.
The upper bound remains untouched but in 2020
Dumitrescu and Mandal improved the lower bound constant to $0.2083$. 
Their approach utilizes the known values of $B_n$ for up to $n=12$.

We tackle the lower bound by utilizing
dynamic programming and the Lindström–Gessel–Viennot lemma. 
Our new bound is $b_n \geq \cOptTwelve$ for sufficiently large~$n$. 
The result is based on a delicate interplay of theoretical ideas and computer assistance.
\end{abstract}

\section{Introduction}

Levi~\cite{Levi1926} introduced arrangements of pseudolines as a natural generalization of line arrangements in 1926. Ringel studied them in the 1950's and Grünbaum~\cite{Gruenbaum1972} popularized them in his 1972 monograph \emph{Arrangements and Spreads}. In the 1980's Goodman and Pollack 
initiated a thorough study and related them to many other objects of interest in discrete geometry.
For a more detailed account to the history and the relevant references we refer the interested reader to the handbook article~\cite{FelsnerGoodman2018}.

An \emph{arrangement of pseudolines} in the Euclidean plane $\R^2$ is a finite family of simple curves, called pseudolines, 
such that each curve approaches infinity in both directions  
and 
every pair intersects in exactly one point 
where the two curves cross. 
More generally, we call a collection of pseudolines \emph{partial arrangement} 
if every pair intersects in at most one crossing-point. 
Pseudolines which do not intersect are said to be \emph{parallel}.
Note that, while for partial arrangements of \emph{proper} lines the relation 'parallel' is transitive, this is no longer true in partial pseudoline arrangements.

In this article, 
the focus will be on \emph{simple} arrangements, 
that is, no three or more pseudolines intersect in a common point (called \emph{multicrossing}). 
Moreover, we consider all arrangements to be \emph{marked}, that is,
they have a unique marked unbounded cell, 
which is called \emph{north-cell}. 
Two arrangements are \emph{isomorphic} if one can be mapped to the other by an orientation preserving homeomorphism of the plane that also preserves the north-cell.

While it is known that the number $B_n$ of non-isomorphic arrangements of $n$ pseudolines grows as $2^{\Theta(n^2)}$, it remains a challenging problem to bound the multiplicative factor of the leading term of $\log_2 B_n = \Theta(n^2)$.
Determining precise values for small values of~$n$ is a challenging task as well, see \cref{table:small_numbers}.
Our focus will be on finding better estimates on the lower bound constant
$
    \clow := \liminf_{n \to \infty} \frac{\log_2 B_n}{n^2} 
$.
One can analogously 
define the upper bound constant $\cupp := \limsup_{n \to \infty} \frac{\log_2 B_n}{n^2}$.
It seems to be open whether $\cupp$ and $\clow$ coincide, i.e., whether the limit~exists.

\input{table_small_numbers.tex}

A lot of work has been done on finding good estimates for $\clow$ and~$\cupp$. 
In the 1980's
Goodman and Pollak~\cite{GoodmanPollack1983} 
investigated pseudopoint configurations, which are dual to pseudoline arrangements, and established
the lower bound $\clow \geq \frac{1}{8}$.
An alternative and slightly easier construction 
for $\clow \geq \frac{1}{12}$  
was given by Matou\v{s}ek in his 
\emph{Lectures on Discrete Geometry} \cite[Chapter~6]{Matousek2002_book}.
As pointed out in~\cite{DumitrescuMandal2020}
one also obtains $\clow \geq \frac{1}{8}$ 
from Matou\v{s}ek's construction via  recursion. 
Concerning the upper bound,
Edelsbrunner, O'Rourke and Seidel  \cite{EOS1986} showed $\cupp < \infty$ via an algorithm that constructs an arrangement
of $n$ pseudolines in $O(n^2)$ time.
It is worth noting that an upper bound
can also be obtained via planar graphs:
The dual graph of arrangement of $n$ pseudolines is planar quadrangulation on $m= \binom{n}{2}+n+1$ vertices. 
Since the number of $m$-vertex planar graphs 
is at most~$2^{O(m)}$ \cite{BGH2003,BGH2011}, 
it follows that $B_n \le 2^{O(n^2)}$.

In the 1990's
Knuth \cite[Section~9]{Knuth1992} improved the bounds to $\clow \geq \frac{1}{6}$ and $\cupp < 0.7925$, and conjectured that $\cupp \le 0.5$.
The upper bound was lowered to $\cupp < 0.6974$ by Felsner~\cite{Felsner1997}.
In 2011,
Felsner and Valtr \cite{FelsnerValtr2011} 
further narrowed the gap by showing $\clow > 0.1887$ and $\cupp < 0.6571$,
and, in 2020 Dumitrescu and Mandal \cite{DumitrescuMandal2020} proved 
the currently best lower bound~$\clow > 0.2083$.

\medskip

In this article, 
we make a substantial step on the lower bound
by proving $\clow > \cOptTwelve$. 

\begin{theorem}
\label{thm:main}
    The number $B_n$ of non-isomorphic simple arrangements of $n$ pseudolines satisfies the inequality $B_n \geq 2^{c n^2 - O(n\log n )}$ with $c > \cOptTwelve$. 
\end{theorem}

\section{Outline}
\label{sec:outline}

Our approach is in the spirit of several previous 
bounds. We consider a specific partial arrangement~$\cal L$ of $n$ lines
consisting of  
$k$ \emph{bundles} $\calL_1,\ldots,\calL_k$ of equally-spaced parallel lines. 
We then define a class of local perturbations to $\cal L$ and consider the number of arrangements that can be obtained by these
perturbations. 
This number is a lower bound 
on $B_n$, and it can be improved by recursively applying the same construction to each 
of the parallel classes~$\calL_i$.

The main difference between the approaches 
lies in the number of bundles $k$ and the choice of locality. 
Matou\v{s}ek and also Felsner and Valtr used three bundles but the locality was increased 
from considering just a triple intersection with its two simple resolutions to the full 
intersection pattern of three bundles.
Dumitrescu and Mandal \cite{DumitrescuMandal2020} increased the number~$k$ of bundles to 12 but restricted the locality to the possible resolutions of the multicrossings. 

\begin{figure}[htb]
    \centering
    \includegraphics{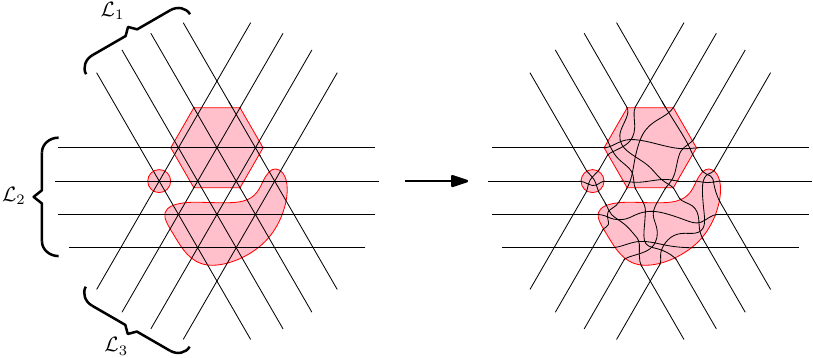}
    \caption{Left: An arrangement of $k=3$ bundles of parallel lines and a collection of interior-disjoint patches (highlighted red) such that each multicrossing point is covered by a patch. 
    Right: A partial pseudoline arrangement with the same parallel bundles obtained by rerouting within the patches.
    }
    \label{fig:phase1}
\end{figure}

Our approach combines higher values of $k$ 
(we discuss cases with $k$ equal to 4, 6, and 12 in this paper) with an increased locality for the perturbations, i.e.,  
we allow reroutings of sets of pseudolines in designated regions, which we call \emph{patches}.
\cref{fig:phase1} gives an illustration.

When rerouting a partial arrangement within a patch~$P$, 
the order of the crossings along the pseudolines may change.
In \cref{sec:dynamic_program} we show that the boundary information of $P$  
determines which pairs of pseudolines cross within~$P$, but the order of 
crossings along the pseudolines is not determined in general.
Outside of~$P$, the arrangement remains unaffected,
which allows us to
count the number of reroutings for each patch independently.
The total number of perturbations is obtained as the product 
of the numbers computed for the individual patches.
Details on how we computed the number of possibilities within a 
patch are given in \cref{sec:rerouting} .

To eventually use computer assistance,
we choose patches of high regularity
and reasonably small complexity.
In fact, since the initial arrangement is extremely regular, it is sufficient 
to determine the rerouting possibilities 
only for a small number of patch-types.
Only a negligible fraction of patches along the boundaries are different.
As we only want to find an asymptotic lower bound on~$B_n$, 
the small number of irregular patches 
along the boundary of the regions 
will be ignored for the counting.

\medskip
\noindent
To eventually prove \cref{thm:main},
we perform the following two steps:

\begin{itemize}\setlength\itemsep{0.5em}

\item 
In the first step (\cref{sec:step1_patches})
we specify the parameters of the construction:
For $k=4$ (\cref{ssec:4slopes}), $k=6$ (\cref{ssec:6slopes}),
and $k=12$ (\cref{ssec:rect12slopes})
we construct the initial arrangement with $k$ bundles of $\floorfrac{n}{k}$ parallel lines
and cover the multicrossing points by patches.
By resolving the multicrossing points within the 
patches, and taking the product over all patches we obtain
a lower bound on the number $F_k(n)$ of partial arrangements with $k$ bundles of $\floorfrac{n}{k}$ parallel pseudolines.

\item
In the second step (\cref{sec:step2_parallelities}),
we recursively account for crossings in bundles of pseudolines 
which had been parallel before.
The product of the so-computed possibilities yields the
improved lower bound on the number $B_n$ of simple arrangements on $n$ pseudolines.
\end{itemize}

\section{Step 1: bundles of parallel lines, patches, and perturbations}
\label{sec:step1_patches}

For the start we fix an integer~$k$ and
construct an arrangement $\calL$ of $k$ bundles of $\floorfrac{n}{k}$ parallel lines
as in~\cite{DumitrescuMandal2020}.
If $n$ is not an odd multiple of~$k$, 
the remaining lines are discarded, or not used
in the counting.
We then cover all multicrossing points
by a family of disjoint patches, and
reroute the line segments within the patches 
so that all multicrossing points are resolved
and the arrangement is simple.
Since we use computers it is convenient to
construct patches with a high regularity.
Moreover, for simplicity in the computational part
(cf.\ \cref{sec:rerouting}) 
we ensure that no crossings of $\calL$ are located on the boundary of a patch
and no patch has a touching with a line.

\begin{itemize}\setlength\itemsep{0.5em}

\item 
\textbf{Extremum 1:}
If we use one tiny patch for each of the multicrossings, 
the counting will give the same results as in \cite{DumitrescuMandal2020},
where each multicrossing was rerouted locally in all possible ways for various configurations with $k$ up to~12.

\item 
\textbf{Extremum 2:}
If we choose one gigantic patch containing all crossings of $\calL$, 
then all \emph{partial} arrangements of $n$ pseudolines with the same parallel bundles as in~$\calL$ will be counted, that is $F_k(n)$.
For the case of $k=3$ bundles,
Felsner and Valtr \cite{FelsnerValtr2011}
determined that $ F_3(n) = 2^{c n^2 + o(n^2)}$ where $c = \frac{\log_2(3)}{2} - \frac{2}{3} \approx 0.1258$. 
\end{itemize}

\subsection{Construction with 4 bundles}
\label{ssec:4slopes}

\begin{figure}[htb]
	\centering
    	
        \hbox{}
        \hfill
        \includegraphics[page=1,width=0.47\textwidth]{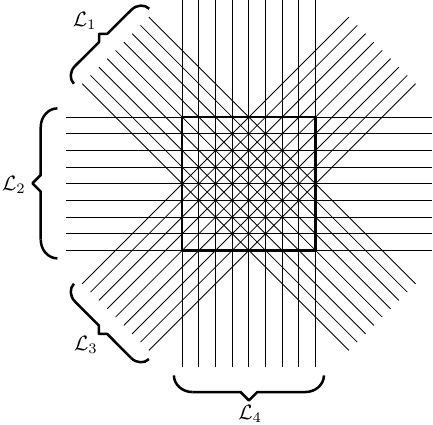}
        \hfill
        \includegraphics[page=2,width=0.47\textwidth]{4slopes_crop.pdf}
        \hfill
        \hbox{}
        
	\caption{Construction with 4 bundles as in \cite{DumitrescuMandal2020}.
 }
	\label{fig:4bundles}
\end{figure}

First we 
construct a partial arrangement $\calL$ of $n=4m$ lines, 
consisting of $4$ bundles 
$\calL_1,\ldots,\calL_4$ 
of~$m$ parallel lines.
See \cref{fig:4bundles} for an illustration.

The construction comes with crossings of order 2, 3, and 4.
We restrict our attention to regions with multicrossings
since regions with 2-crossings do not allow reroutings.
As illustrated on the right-hand side figure,
there are two types of regions with multicrossings:
\begin{itemize}
\item $R_3$ 
contains multicrossings of degree~$3$;
\item $R_4$ 
contains multicrossings of degree~$4$.
\end{itemize}

\begin{figure}[htb]
    \centering
    \includegraphics[page=3,width=0.46\textwidth]{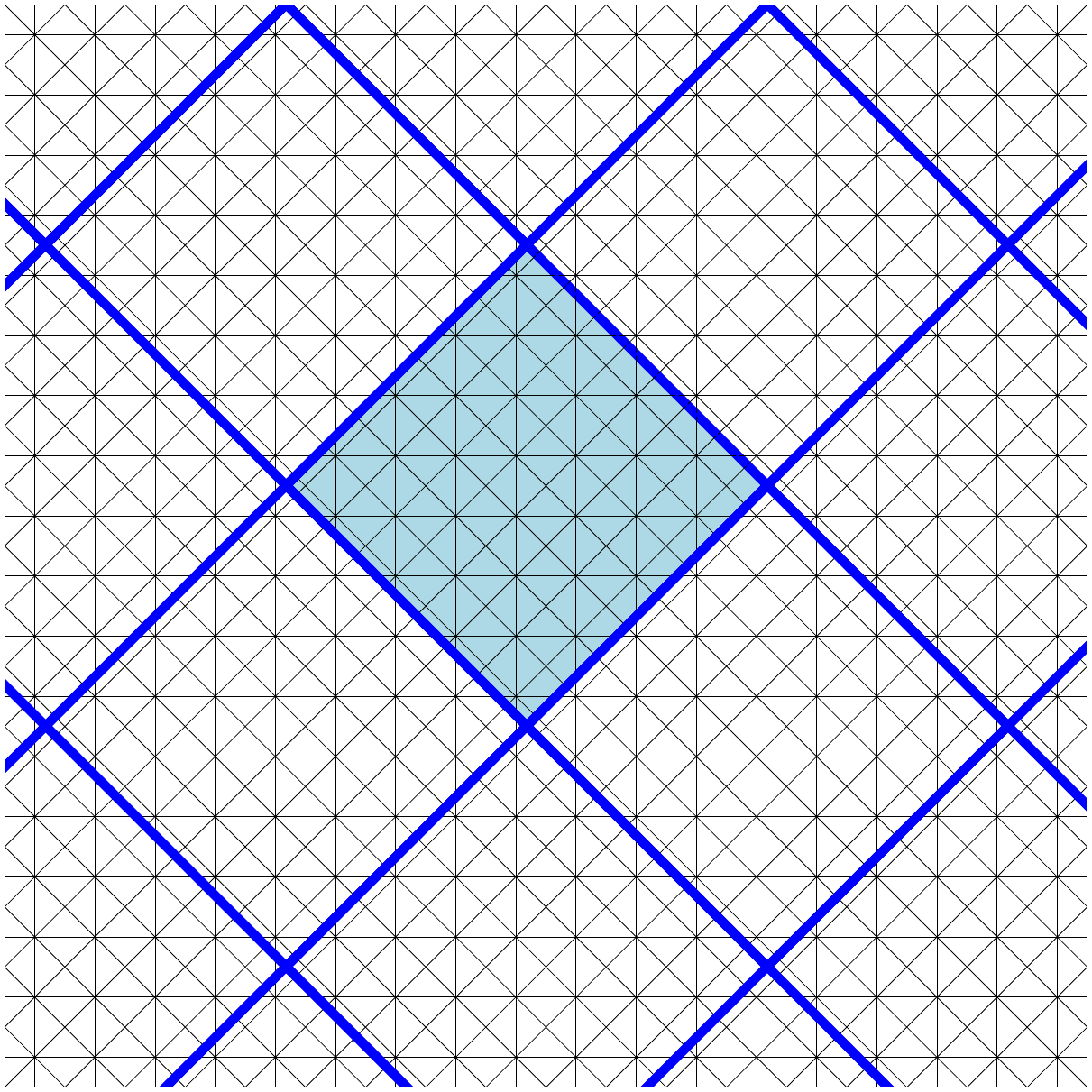}
    
    \caption{
        An illustration of the tiling for region $R_4$ in the 4-slope construction. 
        The patch $P_4$ is a square of side length $8\frac{\sqrt{2}}{2}$ and area~$32$.
    } 
    \label{fig:k4_p4}
\end{figure}

We use different patches in the regions $R_3$ and $R_4$ 
but in both cases rectangular patches can be used to tile the full region apart from a negligible area at its boundary.
See \cref{fig:k4_p4} for the patch $P_4$ we used to cover $R_4$.
Note that each copy of $P_4$ intersects with $\calL$ in an equivalent way.
This will be true implicitly of all patches we consider from now on.
Since $R_3$ has a special pattern of intersection, 
we were able to compute the number of reroutings using a more efficient method 
(see \cref{sec:lgv}) 
than our more general dynamic programming approach 
(see \cref{sec:dynamic_program}). 
This enabled us to use a square of side length $1000$ as our patch $P_3$ for $R_3$, 
which is much larger geometrically compared to $P_4$.

Next we determine 
the numbers~$\mu_i$ of patches of type~$P_i$.
Since the number of crossings in our construction
is asymptotically quadratic in~$n$ and 
each patch contains only a constant number of crossings,
the number of patches of type~$P_i$
is also quadratic.
It is important to note that 
the patches along the boundary of $R_i$ behave differently.
However, since there are only linearly many of these irregular patches, they only affect the lower order error term.
Hence we can omit them in the calculations.

We use a simple area calculation to obtain the $\mu_i$'s. 
First we fix $\area(R_3 \cup R_4) = (m-1)^2$ 
which implies $\area(R_3) = \area(R_4) = \frac{(m-1)^2}{2}$. 
Next we determine the areas of the patches. 
Since at this point all of our patches are rectangular, 
this amounts to determining their side lengths. 
The distance between adjacent lines in $\calL_2$ or $\calL_4$ is $1$ 
and between those in $\calL_1$ or $\calL_3$ it is $\frac{\sqrt{2}}{2}$. 
Thus the areas of the patches are
\begin{itemize}
    \item $\area(P_3) = 1000^2$
    \item $\area(P_4) = (8 \cdot \frac{\sqrt{2}}{2})^2 = 32$.
\end{itemize}
Now we can easily write down the $\mu_i's$:
\begin{itemize}
    \item $\mu_3(P_3,n) = 
    \frac{\area(R_3)}{\area(P_3)} - O(n) 
    = \frac{m^2/2}{1000^2} - O(n) = 
    \frac{n^2}{32000000} - O(n)$
    \item $\mu_4(P_4,n) = \frac{\area(R_4)}{\area(P_4)} - O(n)= 
    \frac{m^2/2}{32}- O(n) = 
    \frac{n^2}{1024} - O(n)$.
\end{itemize}

To compute the number $F(P_4)$ of all possible reroutings within the patches of type $P_4$,
we ran our dynamic program; see \cref{sec:dynamic_program} for details.
To compute the number $F(P_3)$ of all possible reroutings within the patches of type $P_3$,
we made use of the Lindström-Gessel-Viennot lemma;
see \cref{sec:lgv}.
We obtained:
\begin{itemize}
    \item $F(P_3) > 2^{1397192}$
    \item $F(P_4) = 10233480626615962155895931163981261674$.
\end{itemize}
Combining the possibilities from the two types of patches according to the estimate
$F_k(n) \ge \prod_{i=3}^k F(P_{i})^{\mu_i(P_i, n)}$ 
yields 
\begin{proposition}\label{thm:F_4_geq}
    $F_4(n) \ge 2^{cn^2 - O(n)}$ with $c > 0.1637$. 
\end{proposition}

\noindent
By writing  
$
c_i := \lim_{n \to \infty} \frac{\mu_i(n)}{n^2} \cdot \log_2 F(P_i) 
$,
we can see the contributions of the patches $P_3$ and $P_4$ to the leading constant $c = c_3+c_4$ from \cref{thm:F_4_geq}; see \cref{table:summary_k4}.

\input{summary_k4}

\subsection{Construction with 6 bundles}
\label{ssec:6slopes}

\begin{figure}[htb]
	\centering
	
    \includegraphics[width=0.8\textwidth]{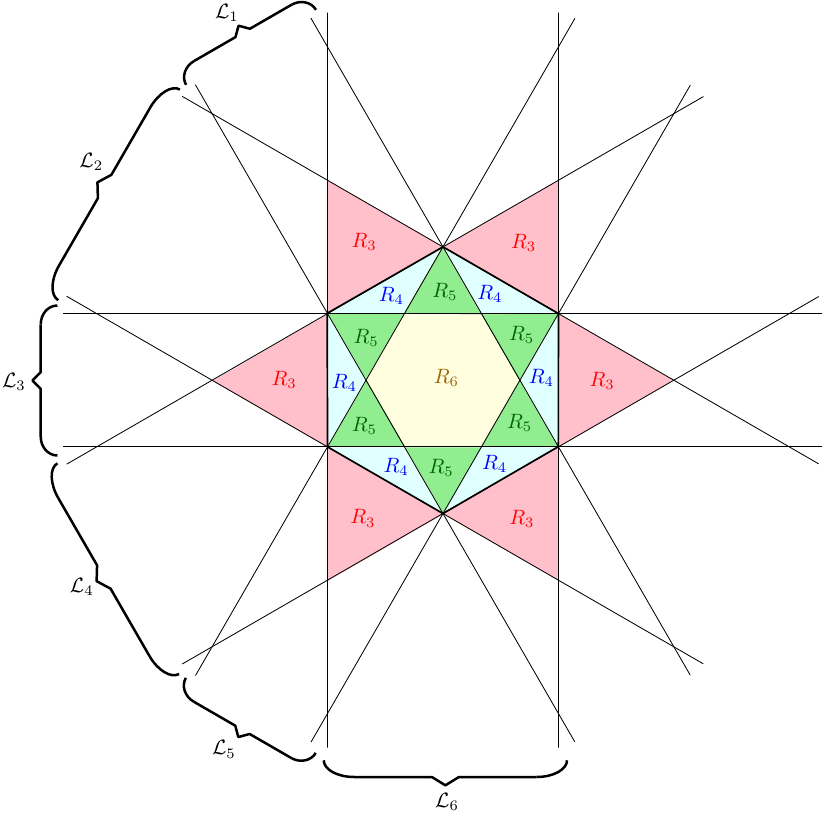}
	
	\caption{Construction with 6 bundles as in \cite{DumitrescuMandal2020}.
 }
	\label{fig:6bundles}
\end{figure}

In this section we consider a partial arrangement $\calL$ of $n$ lines
consisting of  
$6$ bundles 
$\calL_1,\ldots,\calL_6$ 
of $m$ 
parallel lines.
See \cref{fig:6bundles} for an illustration.

\begin{figure}[htb]
    \centering
    
    \subcaptionbox{\label{fig:types_of_tiles_6_slopes_6}}[.31\textwidth]{\includegraphics[width=0.3\textwidth]{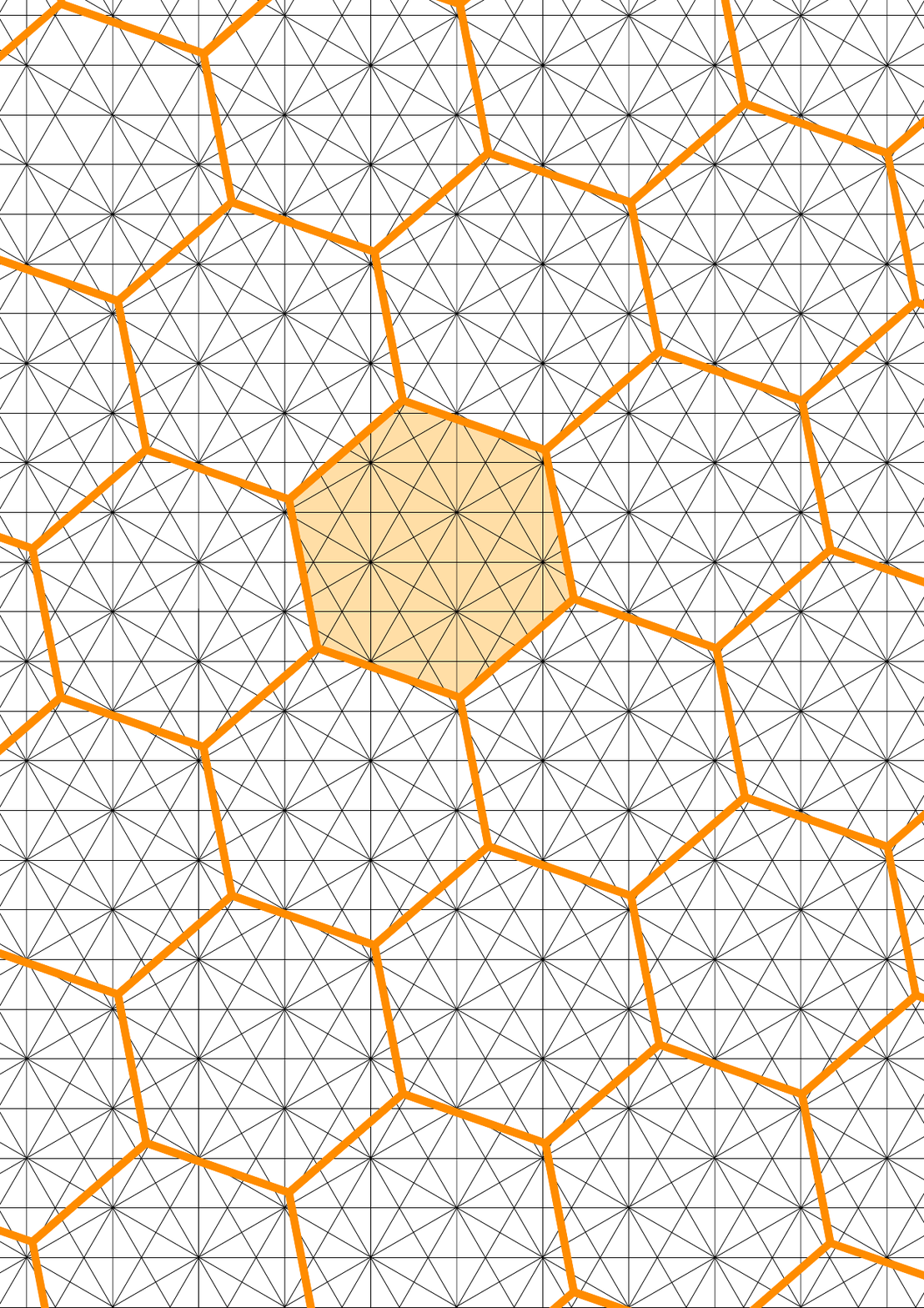}}
    \subcaptionbox{\label{fig:types_of_tiles_6_slopes_5}}[.31\textwidth]{\includegraphics[page=2,width=0.3\textwidth]{hextile_12cr_thick.pdf}}
    \subcaptionbox{\label{fig:types_of_tiles_6_slopes_4}}[.31\textwidth]{\includegraphics[page=3,width=0.3\textwidth]{recttile_thick.pdf}}

    \caption{
    Three of the four types of patches for our construction on $k=6$ bundles:
    \subref{fig:types_of_tiles_6_slopes_6}~For~$R_6$ we use a hexagonal tiling such that each patch $P_6$ contains exactly 7~crossings of order~6 and 14~crossings of order~3.
    \subref{fig:types_of_tiles_6_slopes_5}~For $R_5$ we use a hexagonal tiling such that each patch $P_5$ contains exactly 12~crossings of order~5.
    \subref{fig:types_of_tiles_6_slopes_4}~For $R_4$ we use a rectangular tiling such that each patch $P_4$ contains exactly $8 \cdot 4 = 32$ crossings of order~4.
    } 
    \label{fig:types_of_tiles_6_slopes}
\end{figure}

The construction comes with
four types of regions with multicrossings:
\begin{itemize}
\item 
$R_i$ for $i \in \{3,4,5\}$ 
only contains multicrossings of order~$i$
and 
\item $R_6$ contains multicrossings of order~3 and~6.
\end{itemize}
Note that in contrast to the $k=4$ bundle construction 
from~\cref{ssec:4slopes},
multicrossings of order~3 now occur in $R_3$ and~$R_6$.
For each of the four regions~$R_i$
we will use a different type of patch~$P_i$
that is based on a regular tiling of the plane
to ensure regularity.
\cref{fig:types_of_tiles_6_slopes}
illustrates the tiling for $R_6$, $R_5$ and~$R_4$.
For $R_3$ we use a rhombic tiling such that each patch $P_3$ contains exactly $1000 \cdot 1000$ crossings of order~3; see \cref{sec:lgv}.
Pause to note that
$P_3$ and $P_4$ are affine transformations of the patches used in the 4 slopes
setting in~\cref{ssec:4slopes}.
Therefore they allow the same number of reroutings
$F(P_3)$ and~$F(P_4)$.

We have to determine 
the number $\mu_i$ of patches of type~$i$.
Towards this end we perform a similar area calculation as we did for $k = 4$ bundles. 
This time both our patches and regions are a little more complex and numerous. 
To avoid human error,
we automated the calculation of both area types in different ways.
For the patches, 
we decided to simply count the number of highest degree crossing within the patch,
because that can be done automatically without rounding. 
There is a one-to-one correspondence between such crossings 
and the fundamental parallelogram $\sigma$ of the lattice of crossings between lines from the bundles $\calL_2, \calL_4$ and $\calL_6$.
Since every patch of the same type contains the same number of such crossings
and the patches are part of a potentially infinite tiling,
we obtain the area of $P_i$ 
by simply multiplying the two numbers: 
$\area(P_i) = \#\{\text{crossings of highest degree}\} \cdot \area(\sigma)$.
If we take the area of $\sigma$ to be $1$ 
then the count corresponds exactly to the desired area.
We obtain
\noindent
\begin{align*}
     \area(P_3) = 1000^2, &&
     \area(P_4) = 32, &&
     \area(P_5) = 12, &&
     \area(P_6) = 7.
\end{align*}
Since each region is made up of convex polygons whose supporting lines are known, their areas are easy to compute using computer algebra software. In our supplemental data \cite{supplemental_data}, we provide a Sage script which carries out this task, and can be used for different constructions as well.
The areas are 
\begin{align*}
    \area(R_3) = \tfrac{3(m-1)^2}{4}, &&
    \area(R_4) = \tfrac{(m-1)^2}{4}, &&
    \area(R_5) = \tfrac{(m-1)^2}{4}, &&
    \area(R_6) = \tfrac{(m-1)^2}{4}.
\end{align*}
Using $n=6m$ we can calculate the $\mu_i$'s as we did previously:
\begin{itemize}
    \item $\mu_3(P_3,n) = \frac{\area(R_3)}{\area(P_3)} - O(n) = \frac{3n^2}{144 \cdot 1000^2} - O(n)$
    \item $\mu_4(P_4,n) = \frac{\area(R_4)}{\area(P_4)} - O(n)= \frac{n^2}{144 \cdot 32} - O(n)$
    \item $\mu_3(P_5,n) = \frac{\area(R_5)}{\area(P_5)} - O(n) = \frac{n^2}{144 \cdot 12} - O(n)$
    \item $\mu_3(P_6,n) = \frac{\area(R_6)}{\area(P_6)} - O(n) = \frac{n^2}{144 \cdot 7} - O(n)$.
\end{itemize}

To compute the numbers $F(P_i)$ of all possible perturbations within the patch type $P_i$ for $i=4,5,6$,
we used our dynamic programming approach;
see \cref{sec:dynamic_program} for more details.
For region $F(P_3)$ we were able to use the same number already computed for \cref{ssec:4slopes}. The computation was based on the LGV lemma; see \cref{sec:lgv}. 
Altogether we computed 
\begin{itemize}
    \item $F(P_3) > 2^{1397192}$ 
    \item $F(P_4) = 10233480626615962155895931163981261674$
    \item $F(P_5) = 32207077855497546508132740267$ 
    \item $F(P_6) = 5489259325039519956333073658900$.
\end{itemize}
From $F_k(n) \ge \prod_{i=3}^k F(P_{i})^{\mu_i(n)}$,
we can now derive:
\begin{proposition}\label{thm:F_6_geq}
$F_6(n) \ge 2^{cn^2 - O(n)}$
with $c > 0.2118$.
\end{proposition}

\noindent
By writing  
$
c_i := \lim_{n \to \infty} \frac{\mu_i(n)}{n^2} \cdot \log_2 F(P_i)
$,
we can see the contributions of the patches $P_3$, $P_4$, $P_5$ and $P_6$ to the 
leading constant  $c = c_3+c_4+c_5+c_6$  from \cref{thm:F_6_geq}; see \cref{table:summary_k6}.

\input{summary_k6}

\subsection{Rectangular construction with 12 bundles}
\label{ssec:rect12slopes}

\begin{figure}[t]
    \centering
    \includegraphics[width=0.8\linewidth]{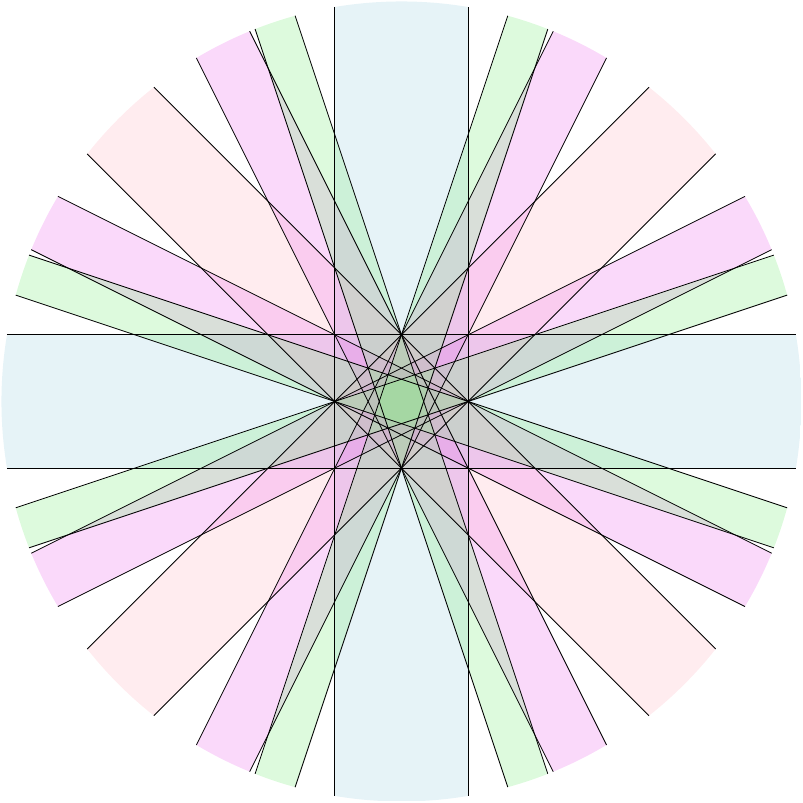}
    \caption{The construction with 12 bundles as in~\cite{DumitrescuMandal2020}.}
    \label{fig:12bundles}
\end{figure}

\begin{figure}[]
    \centering

    \subcaptionbox{\label{fig:regions12_outer}}[.49\textwidth]{\includegraphics[width=0.46\textwidth]{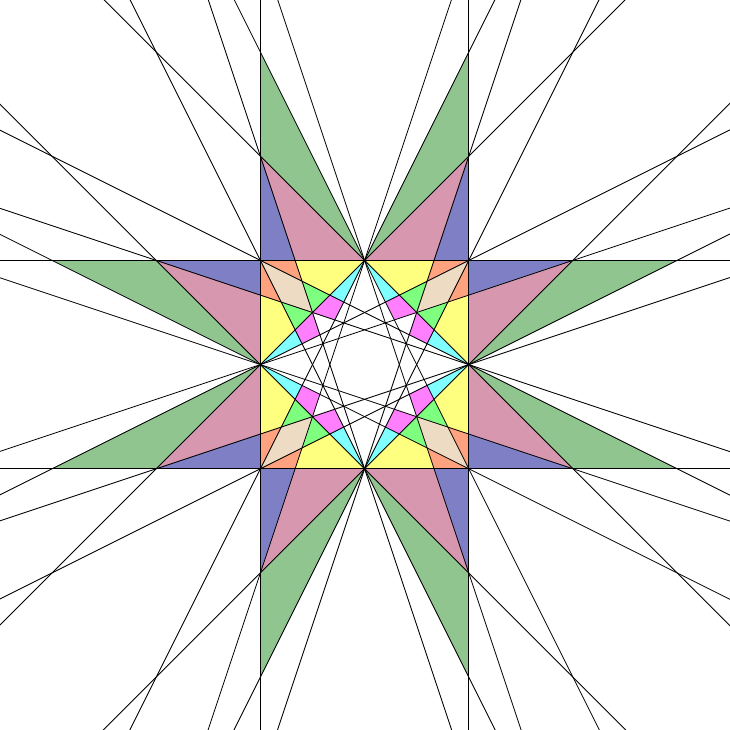}}
    \subcaptionbox{\label{fig:regions12_inner}}[.49\textwidth]{\includegraphics[width=0.46\textwidth]{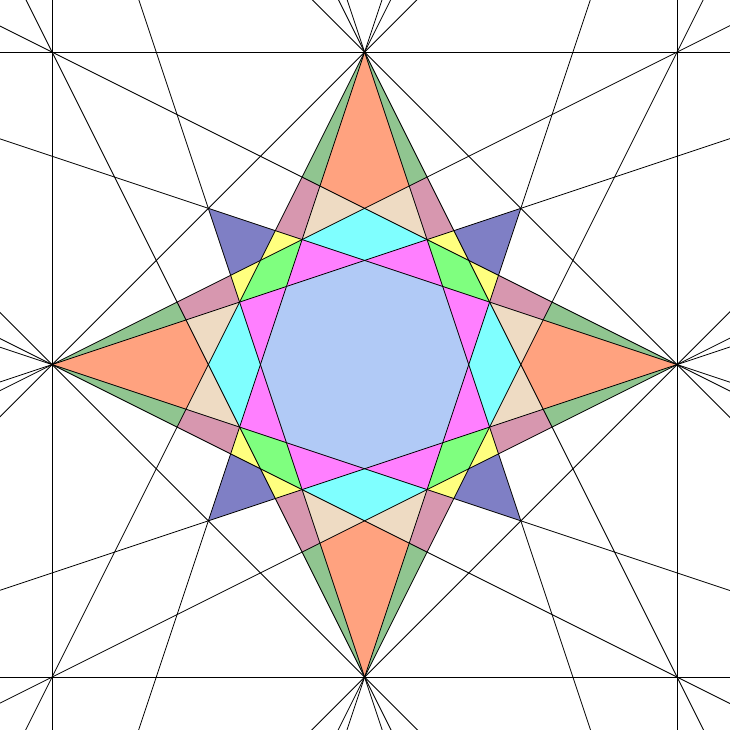}}
        
    \caption{Regions with multi-crossings in the rectangular 12-slope construction. \\
\subref{fig:regions12_outer} shows the nine outer regions $
{\color{color1}R_{3a}}
,
{\color{color2}R_{3b}}
,
{\color{color3}R_{4a}}
,
{\color{color4}R_{4b}}
,
{\color{color5}R_{5a}}
,
{\color{color6}R_{5b}}
,
{\color{color7}R_{6a}}
,
{\color{color8}R_{6b}}
,
{\color{color9}R_{7a}}
$.\\
\subref{fig:regions12_inner} shows the
ten inner regions $
{\color{color1}R_{7b}}
,
{\color{color2}R_{8a}}
,
{\color{color3}R_{8b}}
,
{\color{color4}R_{8c}}
,
{\color{color5}R_{9a}}
,
{\color{color6}R_{9b}}
,
{\color{color7}R_{10a}}
,
{\color{color8}R_{10b}}
,
{\color{color9}R_{11}}
,
{\color{color10}R_{12}}
$.
}
    \label{fig:regions12}
\end{figure}

\begin{figure}[]
    \centering
    \includegraphics[width=0.46\linewidth]{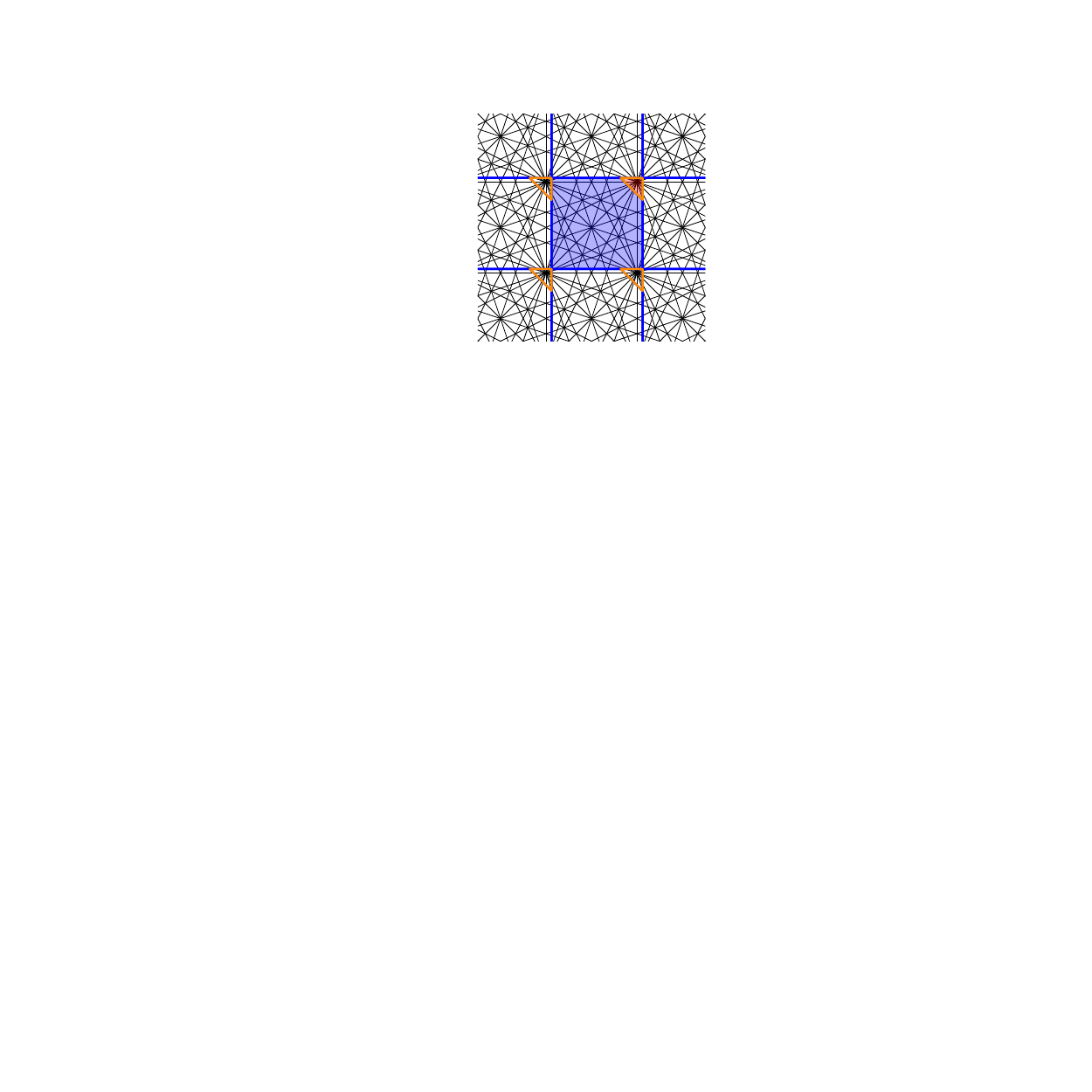}
    \caption{    
    To partition region~$R_{12}$ in the $k=12$ bundle construction, we use a square-shaped patch~$P_{12}$ which is further partitioned into two sub-patches. 
    }
    \label{fig:RA_zoom}
\end{figure}

Finally we consider a partial arrangement $\calL$ of $n$ lines 
consisting of $12$ bundles $\calL_1, \ldots, \calL_{12}$ 
of $m$ parallel lines. 
See \cref{fig:12bundles} for an illustration.

Up to symmetries of $\calL$, 
the construction produces $19$ regions 
that contain intersections of at least three bundles, 
see \cref{fig:regions12}.
Since structurally different regions may meet the same numbers of bundles, 
we differentiate them by adding lowercase Latin letters to their label arbitrarily.
Similar as above, 
we cover each region 
(except for a small area at its boundary)
using a regular tiling 
and count the numbers of reroutings within each tile using a computer. 

For region $R_{12}$ 
even a $1 \times 1$ square
is too complex for our dynamic program to handle in reasonable time/memory.
Therefore, we used two sub-patches as illustrated in \cref{fig:RA_zoom}.
Multiplying the obtained numbers of reroutings in the two sub-patches
yields a lower bound for the $1 \times 1$ square. 
Note that the smaller patch contains only one crossing of degree~$12$. 
Since all pairs of pseudolines cross, 
the number of reroutings is exactly the number $B_{12}$ of arrangements of $12$ pseudolines.

For the other regions we used specific patches which are available in our supplemental data \cite{supplemental_data}. 
To obtain the numbers $\mu_{ix}$ 
of copies of patch $P_{ix}$,
we simply divide the area of the whole region $R_{ix}$ 
by the area of the patch:
\begin{equation*}
    \mu_{ix} = \frac{ \area(R_{ix}) }{ \area(P_{ix}) } - O(n),
\end{equation*}
see \cref{table:summary_k12}.
To compute the areas of the patches we count the number of highest degree crossings they contain again. In this construction, those correspond to $1 \times 1$ squares, so the resulting count will be equal to the patch's area.
The areas of the regions can be computed automatically 
from the specification of the bundles
using a Sage script that we provide with our supplemental data \cite{supplemental_data}. 

Now we compute the numbers $F(P_{ix})$ of possible reroutings within each patch, 
see \cref{table:summary_k12}.
To cope with the large patches in the two regions $R_{3a}$ and $R_{3b}$, 
we again used the Lindström-Gessel-Viennot lemma,
see \cref{sec:lgv}.
For all the other patches we used a dynamic programming approach, 
see \cref{sec:dynamic_program}.
Using the estimate $F_k(n) \ge \prod F(P_{ix})^{\mu_{ix}(n)}$
we derive:
\begin{proposition}\label{thm:F_12_geq}
$F_{12}(n) \ge 2^{cn^2 - O(n)}$
with $c > 0.2494$.
\end{proposition}
By writing  
$c_{ix} := \lim_{n \to \infty} \frac{\mu_{ix}(n)}{n^2} \cdot \log_2 F(P_{ix})$,
we can see the contributions of the patches to the 
leading constant  $c = \sum c_{ix}$  from \cref{thm:F_12_geq}.
All the data discussed above
are summarised in \cref{table:summary_k12}.

\input{summary_k12}

\clearpage

\section{Step 2: resolving parallel bundles}
\label{sec:step2_parallelities}

\begin{figure}[htb]
    \centering
    \includegraphics[width=0.9\textwidth]{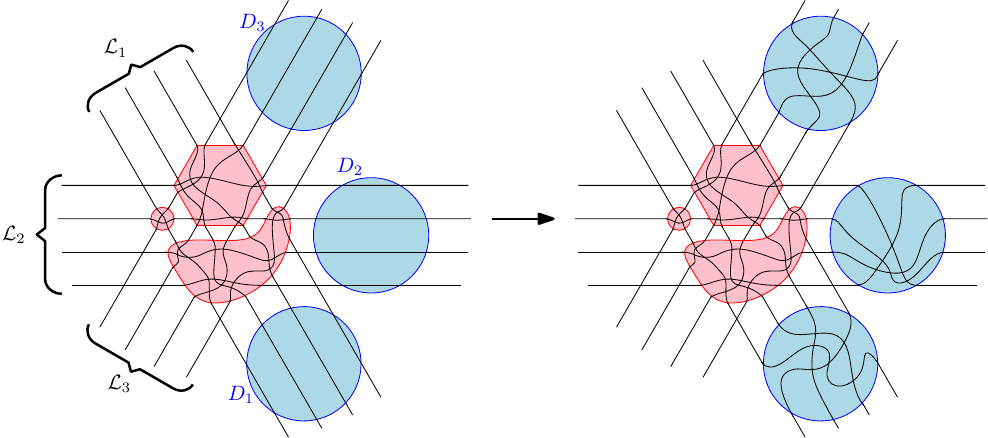}
    \caption{Left: A partial arrangement of $3$ bundles of parallel pseudolines and a collection of interior-disjoint disks (highlighted blue) such that each bundle is covered by one disk. \\
    Right: A proper pseudoline arrangement obtained by rerouting within the disks.
    }
    \label{fig:phase2}
\end{figure}

With the second and final step, 
we obtain a simple arrangement of pairwise intersecting pseudolines from 
a partial arrangement of $k$ bundles of $m = \floorfrac{n}{k}$ parallel pseudolines.
To do so, we use a recursive scheme as in~\cite{FelsnerValtr2011,DumitrescuMandal2020} 
to make each pair of parallel pseudolines cross:
For each $i=1,\ldots,k$, 
we consider a disk $D_i$ such that 
\begin{enumerate}[(1)]
\item $D_i$ intersects all parallel pseudolines of the bundle~$\calL_i$ and no other pseudolines,
and
\item
no two disks overlap.
\end{enumerate}
Within each disk $D_i$ we can place any of the $B_m$ arrangements of $m$ pseudolines. This makes all the pseudolines of a bundle cross.
\cref{fig:phase2} gives an illustration for the case $k=3$.

Since
all $D_i's$ are independent
and there are $B_m$ possibilities to reroute within each $D_i$, we obtain the estimate
\[
B_n \ge \underbrace{F_k(n)}_\text{Step 1} \cdot \underbrace{(B_m)^k}_\text{Step~2},
\]
where $m= \floorfrac{n}{k}$.
With the following lemma
we can derive
$\clow \geq \frac{k}{k-1}c$ 
where $c$ is the constant 
obtained in~\cref{sec:step1_patches}.
The construction with $k=4$ bundles  
yields $\clow > \cOptFour$, 
which is already an improvement 
to the previous best bound by Dumitrescu and Mandal \cite{DumitrescuMandal2020}.
The construction with $k=6$ bundles 
allows a bigger step:
It gives the lower bound $\clow > \cOptSix$.
The biggest improvement comes from the construction with $k=12$ bundles, 
which yields $\clow > \cOptTwelve$
and therefore completes the proof of~\cref{thm:main}.

\begin{lemma}\label{lem:inductive_improvement}
If $F_k(n) \geq 2^{c n^2 - O(n)}$ for some $c > 0$ then 
$B_n \geq 2^{\frac{k}{k-1}c n^2 - O(n\log n)}$. 
\end{lemma}   

The proof of \cref{lem:inductive_improvement} is defered to
\ifthenelse{\boolean{appendix}}
{\cref{appendix}.}
{the full version~\cite{fullversion}.}

\section{Counting all possible reroutings in a patch}
\label{sec:rerouting}

In the following we present two approaches how to compute the number $F(P)$ of all possible perturbations within a patch~$P$.
In \cref{sec:dynamic_program} we discuss a general approach based on dynamic programming which can deal with any patch.
In \cref{sec:lgv} we discuss the special case of 3 slopes, where the reroutings can be derived more efficiently with the Lindström-Gessel-Viennot lemma.

\subsection{Counting reroutings via dynamic programming}
\label{sec:dynamic_program}

\begin{figure}[htb]
	\centering
	
    \includegraphics[page=2,width=0.8\textwidth]{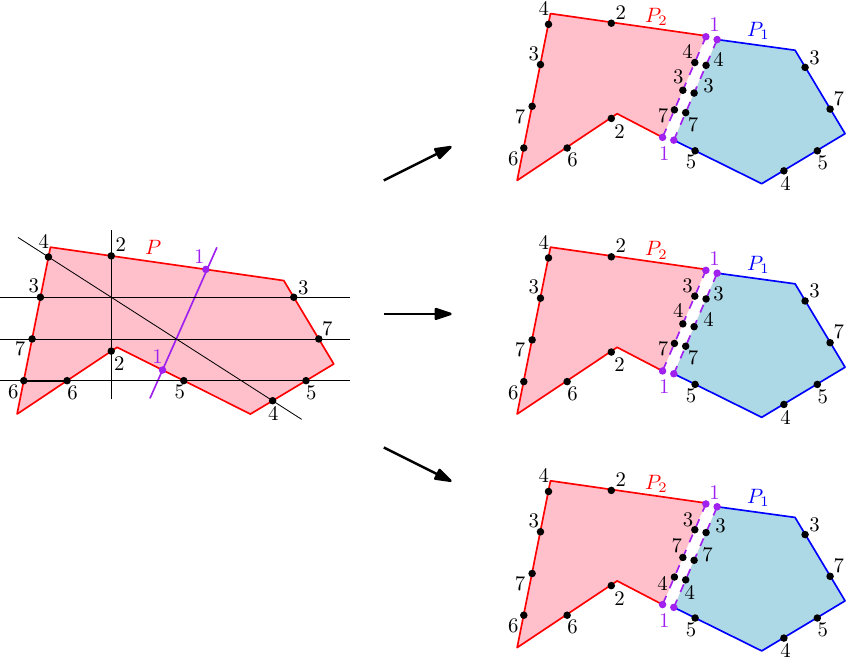}
	
	\caption{
    An illustration of how to recurse on a patch~$P$.
    When cutting along segment~$1$, highlighted purple,
    there are intersections with the segments $3$, $4$, and $7$.
    As the segments $3$ and $7$ do no cross within $P$,
    there are only three possibilities for placing the 
    three crossings along the segment~$1$, namely
    $4$--$3$--$7$ (right top),
    $3$--$4$--$7$ (right center) and
    $3$--$7$--$4$ (right bottom).
 }
	\label{fig:buble_bip}
\end{figure}

Consider a patch~$P$
with boundary curve~$C$.
Since $C$ has no touching with any line,
the intersection of~$P$ with any line $L \in \calL$ 
yields proper line segments.
And because $C$ passes no crossings all the end points of segments are distinct.
We can label the end points on $C$ according to the index of the resulting line segments
and since every index occurs twice, we call the cyclic sequence of end points along~$C$ a \emph{bipermutation}.
Note that in the case of a non-convex boundary curve~$C$, 
the intersection of a line with $P$ consist of  several segments, they have to be distingushed by distinct labels.
\cref{fig:buble_bip} gives an illustration.

The number of possible reroutings within a patch~$P$ is computed recursively:
Choose a segment $z$, it splits $P$ into two parts. For every segment  $a \neq z$ we can 
determine whether $a$ crosses $z$ by looking at the occurrences of $a$ and $z$ 
in the bipermutation. If up to cyclic shifts the pattern is $azaz$
they cross in $P$, if the pattern is $aazz$ 
they are parallel in $P$. 

Now for every legal order~$\pi$ of the crossing segments on $z$ we 
combine the two arcs of the bipermutation defined by the two occurrences 
of $z$ and combine them with $\pi$ to obtain the bipermutations of the two 
smaller patches $P_1$ and~$P_2$.
\cref{fig:buble_bip} gives an example.

The legal orderings of crossings on $z$ can be 
found as follows:
Consider a pair of segments $a,b$ that both cross~$z$. 
Either $a$ and $b$ are parallel in~$P$
(i.e., we see $abzbaz$ in the bipermutation)
or they form a crossing in~$P$ which we denote by $\crossing{a}{b}$
(i.e., we see $abzabz$).
In \cref{fig:buble_bip} 
the segments $3$ and $7$ are parallel,  
and segments $3$ and $4$ cross,
and the crossing point $\crossing{3}{4}$ can lie on either side of the segment~$1$.
In the case where $a$ and $b$ are parallel,
the order of the two crossings $\crossing{a}{z}$ and $\crossing{b}{z}$ along $z$ is uniquely determined.
In the case where $a$ and $b$ cross,
the crossing point $\crossing{a}{b}$
can lie on either side of~$z$,
and the choice of the part further 
determines the order of the two crossings $\crossing{a}{z}$ and $\crossing{b}{z}$ along~$z$.
Altogether, the bipermutation of~$P$ 
determines a partial order $\prec$ of all crossings on~$z$ such that the legal orderings correspond to linear extensions of this partial ordering.

To count all possible reroutings within~$P$,
we iterate over all linear extensions $\prec^\star$ of $\prec$ and continue recursively.
Note that linear extensions can be enumerated efficiently via back-tracking: 
pick a maximal element and recurse on the remaining elements.
Each extension $\prec^\star$ uniquely determines the order of crossings along~$z$ and therefore the two sides $P_1(z,\prec^\star)$ and $P_1(z,\prec^\star)$ can be explicitely given via their bipermutations.
By cumulatively summing up the obtained numbers,
we get
\[
F(P) = \sum_{\text{$\prec^\star$ linear extension of $\prec$}} F(P_1(z,\prec^\star)) \cdot F(P_2(z,\prec^\star)).
\]

We provide a computer-assisted framework \cite{supplemental_data}
that allows to fully automatically 
compute $F(P)$ for a given patch~$P$.
The input is given as IPE-file\footnote{IPE \cite{IPE} is a drawing editor for creating vector graphics in XML or PDF format. Besides the supplemental input files, also all figures in this article were created with it.}. 
The program reads the collection of lines and the polygonal boundary of~$P$,
computes the bipermutation,
and then performs a dynamic program
to determine the number of reroutings within~$P$.
More specifically,
we compute for each bipermutation the lexicographic minimal among all relabelings of the elements,
and reuse previously computed values whenever possible.

\subsection{Counting reroutings via the Lindström-Gessel-Viennot lemma}
\label{sec:lgv}

\begin{figure}
    \centering
    \includegraphics[width=0.75\linewidth]{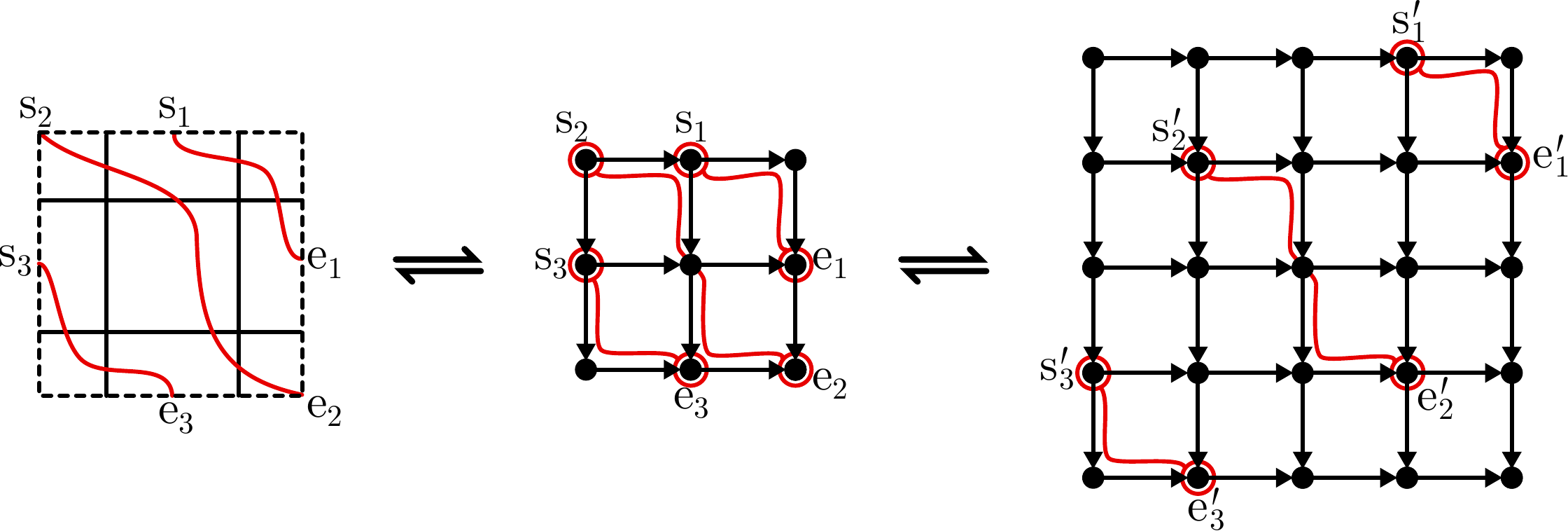}
    \caption{Illustration of the correspondence between arrangements of the chord diagram and disjoint paths in the directed acyclic graph.}
    \label{fig:gessel-viennot}
\end{figure}

We illustrate our approach on a small patch defined by a square of side length $l=2$ (see \cref{fig:gessel-viennot}). We can count the number of reroutings in this patch by starting with only the black grid, and counting the number of combinatorially distinct ways to insert three curves starting at $s_1$, $s_2$ and $s_3$ respectively, and ending at $e_1$, $e_2$ and $e_3$ respectively, such that these three curves do not intersect, and each curve crosses every black segment at most once. We can in turn translate this into counting the number of ways to have three paths in the dual grid graph (with edges oriented left to right and top to bottom), starting and ending at prescribed vertices, such that these paths do not cross (middle illustration in \cref{fig:gessel-viennot}). Note that here by ``not crossing'' we do not mean vertex- or edge-disjoint, but that if a path starts above/right of another, it can at no point go below/left of it. Given three such paths, ordered from top right to bottom left, we can shift the second path one unit down and left along the grid and the second path two units down and left along the grid.
This results in three vertex-disjoint paths (right illustration in \cref{fig:gessel-viennot}). The reverse is also true: given three vertex-disjoint paths on the directed acyclic graph starting at ${s'_1, s'_2, s'_3}$ and ending at ${e'_1, e'_2, e'_3}$, reversing the shifts produces three non-crossing paths starting at ${s_1, s_2, s_3}$ and ending at ${e_1, e_2, e_3}$. 

Thus, our question reduces to counting the number of ways to have three such vertex-disjoint paths, whose set of starting vertices is ${s'_1, s'_2, s'_3}$ and set of ending vertices is ${e'_1, e'_2, e'_3}$.
Note that the vertex-disjointness ensures that in such a collection of paths, 
the path starting at $s'_i$ will end at $e'_i$ for all $1\leq i \leq 3$. 
The Lindström–Gessel–Viennot lemma
gives us an efficient method to do so:
\begin{lemma}[Lindström \cite{Linstrom1973}, Gessel \& Viennot \cite{Gessel1985}]
    Let $G$ be a finite directed acyclic graph. Consider starting vertices $S=\{s_{1},\ldots ,s_{k}\}$ and destination vertices $E=\{e_{1},\ldots ,e_{k}\}$. For any two vertices $u$ and $v$, let $p(u,v)$ be the number of paths from $u$ to $v$. Assume that for any tuple of $k$ vertex-disjoint paths starting in $S$ and ending in $E$, the path starting at $s_i$ necessarily ends at $e_i$, for all $1\leq i \leq k$. Then the number of distinct such tuples is the determinant of the matrix
    \[M =\begin{pmatrix} 
            p(s_1,e_1) & p(s_1, e_2) & \ldots & p(s_1, e_k)\\
            p(s_2,e_1) & p(s_2, e_2) & \ldots & p(s_2, e_k)\\
            \vdots & \vdots & \vdots & \vdots\\
            p(s_k,e_1) & p(s_k, e_2) & \ldots & p(s_k, e_k)\\
        \end{pmatrix}.\]
\end{lemma}

In our specific case, the entries of this matrix are the number of paths going only down or right between two specified vertices in a grid, which can be easily expressed as a binomial coefficient. We have
\begingroup
\renewcommand*{\arraystretch}{1.2} 
\[M =   \begin{pmatrix}
            \binom{2}{1} & \binom{3}{0} &  0\\
            \binom{3}{3} & \binom{4}{2} & \binom{3}{0}\\
            0 & \binom{3}{3} & \binom{2}{1}\\
        \end{pmatrix}
    =   \begin{pmatrix}
            2 & 1 &  0\\
            1 & 6 & 1\\
            0 & 1 & 2\\
        \end{pmatrix}.
        \]
\endgroup

\begin{figure}
    \centering
    \includegraphics[width=0.75\linewidth]{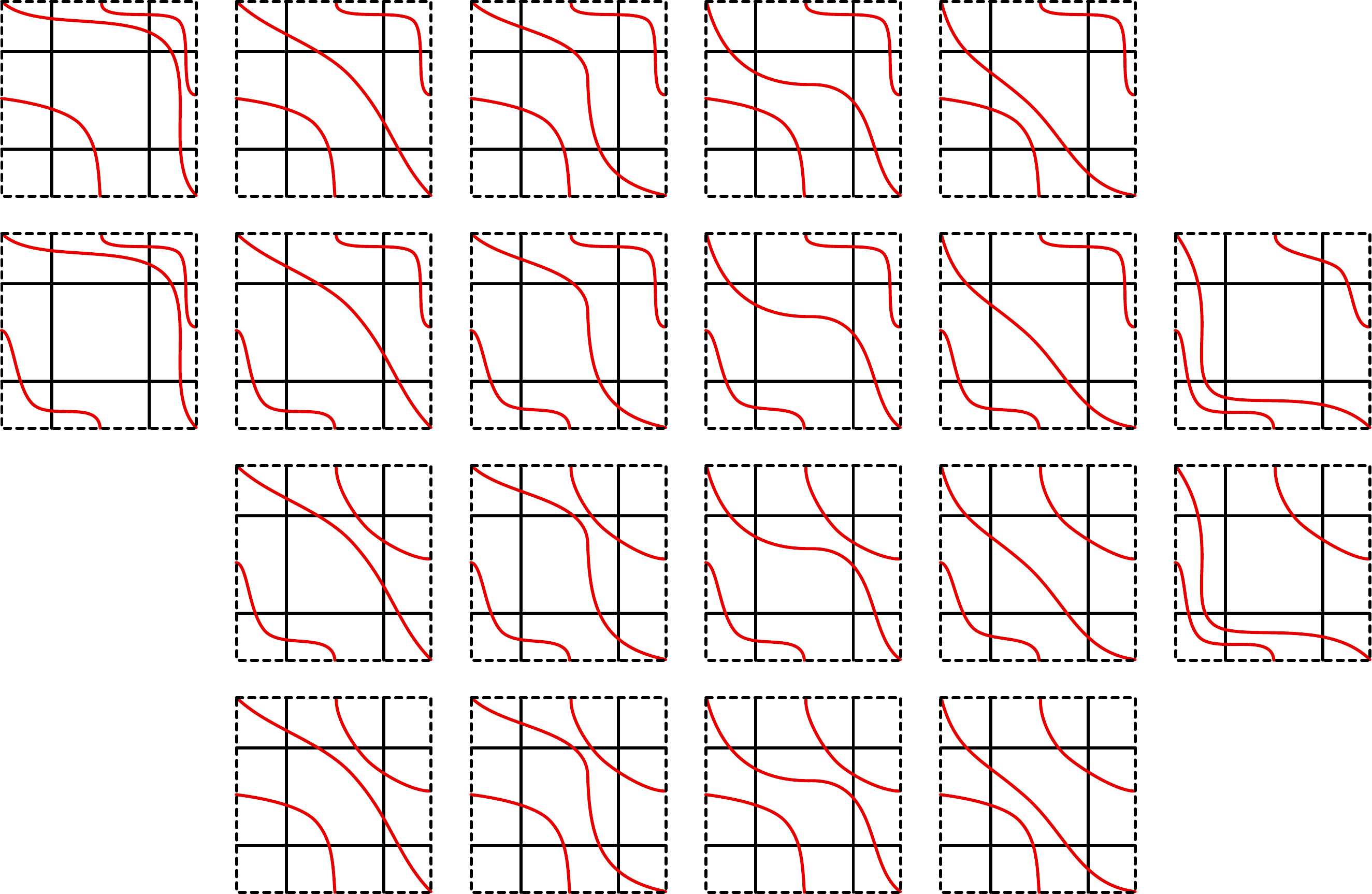}
    \caption{The $20$ reroutings of the patch.}
    \label{fig:LGV-reroutings}
\end{figure}

The determinant of this matrix is $\det(M)=20$, 
corresponding to the $20$ reroutings of the patch illustrated in \cref{fig:LGV-reroutings}.
In general, for a square patch of size $l \times l$, 
the corresponding matrix is $M = (m_{ij})_{1\leq i,j \leq 2l-1}$, with
\[m_{ij} = \binom{2l-|l-i|-|l-j|}{\frac{2l-|l-i|-|l-j|+3|i-j|}{2}}.\]
For our actual bounds, we use a square patch of side length $l=1000$. 
The number of reroutings of this patch is too big to write explicitly here.
In base $10$, it is $420\ 597$ digits long. 
However, it is enough for us to know that $\log_2(\text{\# of reroutings}) > 1397192$. 
\cref{table:table_lgv} shows the numbers and computation times for other side lengths. 
The computations suggest that the ratio $\frac{\log_2(\text{\# of reroutings})}{l^2}$ is converging as $l$ increases. 
Hence, we expect that increasing $l$ only gives negligible improvements.

\input{table_lgv}

\section{Discussion}
\label{sec:discussion}

All results presented in \cref{table:summary_k4,table:summary_k6,table:summary_k12,table:table_lgv} 
were computed
on cluster nodes of TU Berlin with up to 1TB of RAM. The programs can also be run with fewer resources but the dynamic program comes with a memory--time trade-off.
Our programs and further details are available in our supplemental data \cite{supplemental_data}.
We also provide simpler patches 
for which the program only needs few CPU seconds and low RAM.
The simpler patches have also been verified by an independent implementation \cite{DallantVersion}.
They, however, give a slightly worse bound.

We performed quite some experiments to optimize the set of parameters. 
For the 4-slope construction in \cref{ssec:4slopes},
we so far manged to obtain 
$\clow > \cOptFour$
and
for the 6-slope construction in \cref{ssec:6slopes},
$\clow > \cOptSix$.
For the rectangular 12-slope construction in \cref{ssec:rect12slopes}, 
we obtained
$\clow > \cOptTwelve$,
which is the currently best lower bound.
It is interesting that already the 4-slope construction 
gives a significant improvement 
to the previous best bound by Dumitrescu and Mandal \cite{DumitrescuMandal2020}.
In the future we plan to investigate also the construction with 16 slopes, from which we expect an improvement to~0.28.

While the results from \cite{DumitrescuMandal2020}
suggest that 
larger values of~$k$ give better bounds,
the computations get more and more complex.
In fact, as the number $k$ increases,
the number of different shape types 
and the complexity of the patches increases.
Since our program can only deal with patches containing about 30 to 40 segments in reasonable time, 
depending on the structure of crossings within it,
there is a trade-off between 
the number of crossings  
and the number of bundles involved in a patch.

As long as one fixes~$k$,
the counting approach is implicitly limited by~$F_k(n)$, 
which is much smaller than~$B_n$.
Since $F_3(n) = 2^{c n^2 + o(n^2)}$ with $c = \frac{\log_2(3)}{2} - \frac{2}{3} \approx 0.1258$ is known 
\cite{FelsnerValtr2011},
it would be interesting to
determine $\lim_{n \to \infty} \frac{\log_2 F_i(n)}{n^2}$
for $i=4,\ldots,12$.
We wonder
how far from the truth
the constant in \cref{thm:F_4_geq,thm:F_6_geq} is.
However, since our tiling approach asserts a particular global structure on the arrangement, 
it is impossible to approximate $F_k$ -- even under the assumption that we can process arbitrarily large (but finite) patches. 

There is a lot of freedom for choosing the patches for each region,
and we indeed experimented with various shapes.
The ones presented along \cref{sec:step1_patches} 
led to the best results.
As mentioned along \cref{sec:dynamic_program},
when recursively computing~$F(B)$, 
the choice of $z$ can have a significant 
impact on the computing time.
Experiments showed that the best practical performance is obtained by choosing a cutting-segment $z$ such that 
\begin{itemize}
\item the complexity of the larger part $P_2$ is as small as possible.
\end{itemize}
In practice this strategy tends to lead to
\begin{itemize}
\item \emph{balanced} cuts, that is, $P_1$ and $P_2$ are of similar sizes;
\item \emph{short} cuts, that is, the number of intersections on $z$ is relatively small; and
\item relatively small numbers of is legal orders.
\end{itemize}
In general each of those criteria may lead to a different cutting strategy.
However, a properly balanced cut may come with larger parts and hence might be less efficient.
Also the shortest cut will often be very unbalanced, resulting in one half that is only slightly smaller than the whole path and 
thus more recursion steps might be required to deal with the patch.

A good splitting strategy should also take the number of legal orders into account 
because a large number of legal orders will also negatively affect the computation time.
The main advantage of our strategy is that it minimizes the complexity of both $P_1$ and~$P_2$.
This allows the algorithm to reach the bottom of the recursive search tree as fast as possible.
In practice it appeared to be more important to reach a low level in the search than to minimize the number of legal orders
because the number of cache-hits (pre-computed values) increases fast as the level decreases.
To obtain the best possible run-times
we performed benchmarks and made statistics for the computing time used in different layers in the recursive search tree. 
However, in general,
it is hard to tell which cutting strategy is the best because
we do not have any a~priori estimates 
for the computation time or 
the number of reroutings for a patch.

Last but not least,
since local rerouting is a well-known and frequently used technique in combinatorial geometry,
our technique might also be adapted to various other combinatorial structures to derive improved lower bounds.
We see for example great potential to improve the lower bounds for 
arrangements of pseudocircles
\cite{FelsnerScheucher2019},
simple drawings of the complete graph \cite{PachToth2006,Kyncl2009},
or
higher dimensional
pseudohyperplane arrangements \cite[Corollary~7.4.3]{BjoenerLVWSZ1993},
which are all natural generalizations of pseudoline arrangements.

{
	\bibliography{references}
}

\clearpage

\appendix

\ifthenelse{\boolean{appendix}}
{

\section{Proof of \cref{lem:inductive_improvement}}
\label{appendix}

    Choose $L > 0$ as a sufficiently large constant 
    such that $F_k(n) \geq 2^{c n^2 - L n}$ holds for all $n \geq 1$. 
    We show by induction that $B_n \geq 2^{\tilde{c}n^2-G(n)}$ 
    where $G(n) := \tilde{L} n \log_k(\tilde{L} n)$, $\tilde{L} := L + 2\tilde{c}$, and $\tilde{c} := \frac{k}{k-1} \cdot c$.    
    The base case $n \le k$ 
    is clearly satisfied if $L$ is chosen sufficiently large. 
    For the induction step, 
    we have
\begin{align*}
    \log_2 B_n & \geq \log_2\left( F_k(n) \cdot (B_m)^k \right) 
    \\
    &= \log_2  F_k(n)  + k \log_2  B_m 
    \\
            & \stackrel{I.H.}{\geq} (c n^2 - L n) + k \left( \tilde{c} m^2 - G(m) \right)\\
            &\geq c n^2 - L n + \frac{\tilde{c}}{k} (n-k)^2 - kG(m)\\
            &= \underbrace{(c + \frac{\tilde{c}}{k})}_{= \tilde{c}} n^2 - \underbrace{(L + 2\tilde{c})}_{=\tilde{L}} n - k G(m) +  \underbrace{ck}_{\ge 0},
\end{align*}
where $m=\floorfrac{n}{k}$.
Since $kG(m) + \tilde{L}n \le G(n)$ holds by definition,
this completes the proof of \cref{lem:inductive_improvement}.

}{} 

\end{document}

%% file: table_small_numbers.tex
\begin{table}[htb]
\centering

        \begin{tabular}{lrl}
            $n$ & $B_n$\\
            \hline
            3 & 2                               & \\
            4 & 8                               & \\
            5 & 62                              & \\
            6 & 908                             & \\
            7 & 24698                           & \\
            8 & 1232944                         & \\
            9 & 112018190                       & \cite[Section~9]{Knuth1992}\\
            10 & 18410581880                    & \cite{Felsner1997}\\
            11 & 5449192389984                  & \cite{YamanakaNMUN2010}\\
            12 & 2894710651370536               & \cite{Samuel2011}\\
            13 & 2752596959306389652            & \cite{Kawahara2011}\\
            14 & 4675651520558571537540         &  \multirow{2}{*}{Yuma Tanaka 2013}
            \\
            15 & 14163808995580022218786390    \\
            16 & 76413073725772593230461936736  & Günter Rote 2021 \cite{Rote2024}\\
        \end{tabular}
        
\medskip

\caption{The known values of~$B_n$. 
Sequence \href{https://oeis.org/A006245}{A6245} in the OEIS \cite{oeis}.
}

\label{table:small_numbers}

\end{table}

%% file: summary_k4.tex
\begin{table}[htb]
\centering
\begin{tabular}{ l r r r r r}
region  &$\log_2(\text{\# of reroutings})$     &\# of patches     &contribution  &computing time    \\ 
\hline
     $R_3$*& 1397192.00& 1/32000000&       0.04366&      -\\
      $R_4$&     122.94&     1/1024&       0.12006&  20199.30s\\
\hline 
       $\sum$&          -&            -&       0.16373&  20199.30s\\
\end{tabular}

\medskip

\caption{Summary of our computational results for each of the patches in the $k=4$ slope construction.
The first entry (marked with a star) 
is computed using the LGV lemma;
see \cref{table:table_lgv} for details.
The first column shows the name of the region;
the second column shows the number of reroutings ($\log_2(F)$);
the third column shows the number of patches ($\lim \frac{\mu}{n^2}$);
the forth column shows the contribution of the region 
to the leading constant, i.e., the product of the second and third column;
the last column shows the computing time in CPU seconds.
}

\label{table:summary_k4}

\end{table}

%% file: summary_k6.tex
\begin{table}[htb]
\centering
\begin{tabular}{ l r r r r r}
region  &$\log_2(\text{\# of reroutings})$     &\# of patches     &contribution  &computing time    \\ 
\hline
     $R_3$*& 1397192.00&  1/48000000&       0.02911&       -\\
      $R_4$&     122.94&      1/4608&       0.02668&   18787.90s\\
      $R_5$&      94.70&      1/1728&       0.05480&    4219.49s\\
      $R_6$&     102.11&      1/1008&       0.10130&  868737.00s\\
\hline 
       $\sum$&          -&            -&       0.21190&  891744.39s\\
\end{tabular}

\medskip

\caption{Summary of our computational results for each of the patches in the $k=6$ slope construction.
See \cref{table:summary_k4} for description of columns.
}

\label{table:summary_k6}

\end{table}

%% file: summary_k12.tex
\begin{table}[htb]
\centering

\begin{tabular}{ l r r r r r}
region  &$\log_2(\text{\# of reroutings})$     &\# of patches     &contribution  &computing time    \\ 
\hline
     $R_{3a}$*& 1397192.00&  1/144000000&       0.00970&       -\\
     $R_{3b}$*& 1397192.00&  1/432000000&       0.00323&       -\\
     $R_{4a}$&     122.94&       1/6912&       0.01779&   18345.70s\\
     $R_{4b}$&     122.94&      1/69120&       0.00178&   18345.70s\\
     $R_{5a}$&      94.70&       1/6480&       0.01461&    4463.95s\\
     $R_{5b}$&      94.70&      1/17280&       0.00548&    4463.95s\\
     $R_{6a}$&     111.79&      1/17280&       0.00647&   98184.90s\\
     $R_{6b}$&     113.74&      1/17280&       0.00658&   38082.70s\\
     $R_{7a}$&      93.74&       1/7560&       0.01240&    1204.09s\\
     $R_{7b}$&      91.12&      1/20160&       0.00452&    8392.94s\\
     $R_{8a}$&     148.77&      1/24192&       0.00615&  148199.00s\\
     $R_{8b}$&     109.24&      1/15120&       0.00722&   57856.50s\\
     $R_{8c}$&     105.70&       1/6048&       0.01748&  142085.00s\\
     $R_{9a}$&     107.55&      1/30240&       0.00356&    8577.60s\\
     $R_{9b}$&     102.98&      1/10080&       0.01022&     316.13s\\
     $R_{10a}$&     144.70&      1/17280&       0.00837&   68991.30s\\
     $R_{10b}$&      70.78&       1/4320&       0.01639&    1430.85s\\
     $R_{11}$&      99.03&       1/4320&       0.02292&    2059.67s\\
     $R_{12}$&     128.89&       1/1728&       0.07459&      47.15s\\
     \hline 
    $\sum$&          -&            -&       0.24946&  621047.12s\\
\end{tabular}

\medskip

\caption{Summary of our computational results for each of the patches in the $k=12$ slope construction.
See \cref{table:summary_k4} for description of columns.
}

\label{table:summary_k12}

\end{table}

%% file: table_lgv.tex
\begin{table}[htb]
\centering

\begin{tabular}{ r r r r}
$l$  &$\log_2(\text{\# reroutings})$     &ratio    &computing time    \\ 
\hline
10 &   130.523&  1.3052&   0.02s\\
20 &   539.561&  1.3489&   0.09s\\
50 &   3444.189&  1.3777&   0.62s\\
100&   13877.972&  1.3878&   3.44s\\
200&   55719.146&  1.3930&  37.96s\\
500&   349033.080&  1.3961& 2473.19s\\
1000&  1397192.106&    1.3972&  66467.25s
\end{tabular}

\medskip

\caption{Summary on the LGV computations for $l \times l$ square patches for different sizes~$l$.}

\label{table:table_lgv}

\end{table}